\documentclass[12pt]{article}
\usepackage{graphicx}
\usepackage{graphics}
\usepackage{amsthm}
\usepackage{amssymb, amsmath}
\title{Note on the Stieltjes constants: series with Stirling numbers of the first kind}  
\author{Mark W. Coffey\\
Department of Physics\\
Colorado School of Mines\\
Golden, CO  80401\\
mcoffey@mines.edu\\
(Received $\mbox{~~~~~~~~~~~~~~~~~~~~~~~~~~~~~~~2015}$)}
\date{April 2, 2015}
\begin{document}
\maketitle
\baselineskip=25 pt
\begin{abstract}

The Stieltjes constants $\gamma_k(a)$ appear as the coefficients in the regular
part of the Laurent expansion of the Hurwitz zeta function $\zeta(s,a)$ about $s=1$.
We generalize the integral and Stirling number series results of \cite{blago2015} for
$\gamma_k(a=1)$.  Along the way, we point out another recent asymptotic development
for $\gamma_k(a)$ which provides convenient and accurate results for even modest values of $k$.


\end{abstract}
 
\vspace{.25cm}
\baselineskip=15pt
\centerline{\bf Key words and phrases}
\medskip 

\noindent

Stieltjes constants, Riemann zeta function, Hurwitz zeta function, Laurent expansion,
digamma function, polygamma function, harmonic numbers 

\vfill
\centerline{\bf 2010 AMS codes}   
11M35, 11M06, 11Y60.  secondary:  05A10

\baselineskip=25pt
\pagebreak
\medskip
\centerline{\bf Introduction and statement of results}
\medskip


The Stieltjes (or generalized Euler) constants $\gamma_k(a)$ appear as 
expansion coefficients in the Laurent series for the Hurwitz zeta function 
$\zeta(s,a)$ about its simple pole at $s=1$ 
\cite{briggs,coffeyjmaa,coffeydec11,kluyver,mitrovic,stieltjes,wilton2},
$$\zeta(s,a)={1 \over {s-1}}+\sum_{n=0}^\infty {{(-1)^n} \over {n!}}\gamma_n(a)(s-1)^n.
\eqno(1.1)$$
These constants are important in analytic number theory and elsewhere, where
they appear in various estimations and as a result of asymptotic analyses, being given
by the limit relation
$$\gamma_k(a)=\lim_{N \to \infty} \left[\sum_{j=1}^N {{\ln^k (j+a)} \over j}-{{\ln^{k+1}(N+a)} \over {k+1}}\right].$$
In particular, $\gamma_0(a)=-\psi(a)$, where $\psi(z)=\Gamma'(z)/\Gamma(z)$ is the
digamma function, with $\Gamma(z)$ the Gamma function.
With $\gamma$ the Euler constant and $\gamma_1=\gamma_1(1)$ and $\gamma_2=\gamma_2(1)$, we
recall the connection with sums of reciprocal powers of the nontrivial zeros $\rho$ of the
Riemann zeta function,
$$\sum_\rho {1 \over \rho^2}=1-{\pi^2 \over 8}+2\gamma_1+\gamma^2, ~~~~
\sum_\rho {1 \over \rho^3}=1-{7 \over 8}\zeta(3)+\gamma^3+3\gamma \gamma_1+{3 \over 2}\gamma_2,$$
such relations following from the Hadamard factorization. 
We recall the connection of differences of Stieltjes constants with logarithmic sums,
$$\gamma_\ell(a)-\gamma_\ell(b)=\sum_{n=0}^\infty \left[{{\ln^\ell(n+a)} \over 
{n+a}}-{{\ln^\ell(n+b)} \over {n+b}}\right].  \eqno(1.2)$$

An effective asymptotic expression for $\gamma_k$ \cite{knessl1} and $\gamma_k(a)$
\cite{knessl2} for $k \gg 1$ has recently been given.  From these expressions, 
which show accurate estimations for even modest values of $k$,
previously known results on sign changes within the sequence of Stieltjes constants follow.
The asymptotic expressions encapsulate both the magnitude $|\gamma_m(a)|$ and 
the changes in sign of the sequence $\{\gamma_m(a)\}$.  Additionally in \cite{knessl2}, an 
elaborate analysis is provided for a certain alternating binomial sum of the Stieltjes constants.

Evaluations of the first and second Stieltjes constants at rational
argument have been given very recently \cite{blago,coffeyrama}.  These decompositions are effectively Fourier series, thus implying many
extensions and applications, and they supplement the relations presented in 
\cite{coffeystdiffs}.  
Besides elaborating on a multiplication formula for the Stieltjes constants, \cite{coffeyrama} provides examples of integrals evaluating in
terms of differences of the first and second of these constants.  In addition, presented
there is a novel method of determining log-log integrals with a certain polynomial 
denominator integrand. 


In this note, we briefly describe how the integral and series representations in \cite{blago2015}
for $\gamma_k(a=1)$ may be readily generalized to $\gamma_k(a)$.  Having the parameter $a$ is
very useful, as even from the $\gamma_0(a)$ case representations then follow for $\ln \Gamma(a)$,
the polygamma functions, and, using integer arguments, harmonic and generalized harmonic numbers.
The key starting point integral representation of these developments is given, so that proofs are
largely omitted.

The Hurwitz zeta function, initially defined by $\zeta(s,a)=\sum_{n=0}^\infty
(n+a)^{-s}$ for $\mbox{Re} ~s>1$, has an analytic continuation to the whole
complex plane \cite{berndt,titch}.   
In the case of $a=1$, $\zeta(s,a)$ reduces to the Riemann zeta function
$\zeta(s)$ \cite{edwards,riemann}.
In this instance, by convention, the Stieltjes constants
$\gamma_k(1)$ are simply denoted $\gamma_k$ \cite{briggs,kluyver,kreminski,mitrovic,yue}.  We recall that
$\gamma_k(a+1)=\gamma_k(a)-(\ln^k a)/a$, and more generally that for $n \geq 1$ an
integer 
$$\gamma_k(a+n)=\gamma_k(a)-\sum_{j=0}^{n-1} {{\ln^k(a+j)} \over {a+j}},$$
as follows from the functional equation $\zeta(s,a+n)=\zeta(s,a)-\sum_{j=0}^{n-1}
(a+j)^{-s}$.  
In fact, an interval of length $1/2$ is sufficient to characterize the $\gamma_k(a)$'s \cite{hansen}.

Unless specified otherwise below, letters $j$, $k$, $\ell$, $m$, $n$, and $r$ denote positive integers. The Euler constant is given by $\gamma=-\psi(1)=\gamma_0(1)$.  The polygamma
functions are denoted $\psi^{(n)}(z)$ and we note that $\psi^{(n)}(z)=(-1)^{n+1}n!\zeta(n+1,z)$
\cite{nbs,grad}.  We let $s(n,k)$ denote the Stirling numbers of the first kind.  We note the
cutoff property $s(n,k)=0$ for $n<k$, consistent with the combinatorial interpretation of
$(-1)^{n+k}s(n,k)$ as the number of permutations of $n$ symbols having $k$ cycles. 


Reference \cite{blago2015} does not seem to be aware of the representation \cite{coffeystv2}
$$\gamma_k(a)={1 \over {2a}}\ln^k a-{{\ln^{k-1}a} \over {k+1}}+{1 \over a}\int_0^\infty
{{(y/a-i)\ln^k(a-iy)+(y/a+i)\ln^k(a+iy)} \over {(1+y^2/a^2)(e^{2\pi y}-1)}}dy. \eqno(1.3)$$
From this integral representation the following results follow.

{\bf Proposition 1}.  Let Re $a>0$.  Then
(a)
$$\gamma_0(a)=-\psi(a)={1 \over {2a}}-\ln a+{1 \over {2\pi ia}}\int_0^1\left[{1 \over {1-{{\ln(1-u)}
\over {2\pi i a}}}}-{1 \over {1+{{\ln(1-u)}\over {2\pi i a}}}}\right]{{du} \over u},$$
(b) 
$$\gamma_0(a)={1 \over {2a}}-\ln a-{1 \over {\pi a}}\sum_{k=0}^\infty {{(2k+1)!(-1)^k} \over 
{(2\pi a)^{2k+1}}}\sum_{n=1}^\infty {{s(n,2k+1)(-1)^n} \over {n! n}},$$
%
and (c)
$$\gamma_0(a)={1 \over {2a}}-\ln a-{1 \over {2\pi^2 a^2}}\sum_{n=0}^\infty {{(2n+1)!(-1)^n} \over 
{(2\pi a)^{2n}}}\zeta(2n+2).$$

{\bf Proposition 2}.  Let Re $a>0$.  Then (a)
$$\gamma_1(a)={{\ln a} \over {2a}}-{1 \over 2}\ln^2a+\ln a\left[-\psi(a)-{1 \over {2a}}+\ln a\right]$$
$$-{1 \over {\pi a}}\sum_{n=1}^\infty {1 \over {nn!}}\sum_{k=0}^{[n/2]}{{(-1)^k} \over {(2\pi a)^{2k+1}}}|s(2k+2,2)| |s(n,2k+1)|,$$
wherein $s(2k+2,2)=(2k+1)!H_{2k+1}$, $H_n=\sum_{k=1}^n 1/k$ being the $n$th harmonic number,
and (b)
$$\gamma_m(a)={{\ln^m a} \over {2a}}-{1 \over {m+1}}\ln^{m+1} a$$
$$+{1 \over {2\pi i}}\int_0^1\left[{{\ln^m\left(a-{{\ln(1-u)} \over {2\pi i}} \right)}\over {a-{{\ln(1-u)}\over {2\pi i}}}}-{{\ln^m\left(a+{{\ln(1-u)} \over {2\pi i}} \right)} \over {a+{{\ln(1-u)}\over {2\pi i}}}}\right]{{du} \over u}.$$

As we discuss, Proposition 1 may be considered the expansion of $\gamma_0(a)$ or of the digamma
function about $|a| \to \infty$ with arg $a<\pi$.  

\medskip
\centerline{\bf Discussion}
\medskip

The proofs of Propositions 1 and 2 follow from (1.3) and the approach of \cite{blago2015}.
Rather than provide such details, we give a complementary illustrating verification of 
Proposition 1.  For this purpose we recall the generating function (e.g., \cite{nbs})
$$\sum_{n=0}^\infty s(n,k){x^n \over {n!}}={1 \over k!}\ln^k(1+x), \eqno(2.1)$$
holding for $|x|<1$ and $|x|=1$ but excluding $x=-1$.   

{\it Proof of Proposition 1}.  By geometric series, the right side of (a) is given by
$${1 \over {2a}}-\ln a+{1 \over {2\pi ia}}\int_0^1\left[{1 \over {1-{{\ln(1-u)}
\over {2\pi i a}}}}-{1 \over {1+{{\ln(1-u)}\over {2\pi i a}}}}\right]{{du} \over u}$$
$$={1 \over {2a}}-\ln a+{2 \over {2\pi ia}}\int_0^1 \sum_{\stackrel{j=0}{j ~\rm{odd}}}^\infty {{\ln^j(1-u)} \over {(2\pi i a)^j}}{{du} \over u},$$
the integral being evaluated with a change of variable,
$$\int_0^1 {{\ln^j(1-u)} \over u}du=(-1)^{j+1}\int_0^\infty {v^j \over {e^v-1}}dv=(-1)^{j+1}
j!\zeta(j+1).$$
Now \cite{grad} (p.\ 943)
$$\psi(z)=\ln z-{1 \over {2a}}-2\int_0^\infty {{t dt} \over {(t^2+z^2)(e^{2\pi t}-1)}}, ~~~~
\mbox{Re} ~z>0.$$
Then again by the use of geometric series,
$$\psi(z)=\ln z-{1 \over {2a}}-{1 \over {2\pi^2 z^2}}\sum_{j=0}^\infty (-1)^j \int_0^\infty
{v^{2j+1} \over {(4\pi^2 z^2)^j}} {{dv} \over {(e^v-1)}}$$
$$=\ln z-{1 \over {2a}}-{1 \over {2\pi^2 z^2}}\sum_{j=0}^\infty {{(-1)^j} \over {(4\pi^2 z^2)^j}}
(2j+1)!\zeta(2j+2).$$
Thus parts (a) and (c) follow.

Part (b) follows from geometric series expansion of the right side of (a) and the use of the
generating function (2.1).  In particular, the term with $s(0,j)=\delta_{0j}$, the Kronecker
symbol, does not contribute. \qed

{\it Remarks}.  By reordering the double sum and using the cutoff property of Stirling numbers 
of the first kind, part (b) may be rewritten as
$$\gamma_0(a)={1 \over {2a}}-\ln a-{1 \over {\pi a}}\sum_{n=1}^\infty {1 \over {n! n}}
\sum_{k=0}^{[n/2]} {{(2k+1)!(-1)^k} \over {(2\pi a)^{2k+1}}} s(n,2k+1)(-1)^n.$$

The Stirling numbers of the first kind may be written with the generalized harmonic
numbers, and the first few are given by $s(n+1,1)=(-1)^n n!$, $s(n+1,2)=(-1)^{n+1}n!H_n$,
$s(n+1,3)=(-1)^n {{n!} \over 2}[H_n^2-H_n^{(2)}]$, and $s(n+1,4)=(-1)^{n+1}{{n!} \over 6}
[H_n^3-3H_nH_n^{(2)}+2H_n^{(3)}]$.    
The generalized harmonic numbers in terms of polygamma function values
$$H_n^{(r)}=\sum_{k=1}^n {1 \over k^r} ={{(-1)^{r-1}} \over {(r-1)!}}\left[\psi^{(r-1)}(n+1)
-\psi^{(r-1)}(1)\right]={{(-1)^{r-1}} \over {(r-1)!}}\int_0^1{{(t^n-1)} \over {t-1}}\ln^{r-1}t
~dt$$  
enter the representations of Proposition 2 for the higher Stieltjes constants.



\pagebreak

\end{document}